\documentclass[%
reprint,
superscriptaddress,
showpacs,
 amsmath,amssymb,
 aps,
 pra,
]{revtex4-1}

\usepackage{graphicx}
\usepackage{dcolumn}
\usepackage{bm}
\usepackage[dvips]{color}

\begin{document}

\title{A series solution of the general Heun equation in terms of incomplete Beta functions}
\author{A.M. Manukyan}
\affiliation{Institute for Physical Research, NAS of Armenia, Ashtarak, 0203 Armenia}
\author{T.A. Ishkhanyan}
\affiliation{Institute for Physical Research, NAS of Armenia, Ashtarak, 0203 Armenia}
\affiliation{Moscow Institute of Physics and Technology, Dolgoprudny, 141700 Russia}
\author{M.V. Hakobyan}
\affiliation{Institute for Physical Research, NAS of Armenia, Ashtarak, 0203 Armenia}
\affiliation{Yerevan State University, 1 Alex Manookian, Yerevan, 0025 Armenia}
\author{A.M. Ishkhanyan}
\affiliation{Institute for Physical Research, NAS of Armenia, Ashtarak, 0203 Armenia}


\begin{abstract}
We show that in the particular case when a characteristic exponent of the singularity at infinity is zero the solution of the general Heun equation can be expanded in terms of the incomplete Beta functions. By means of termination of the series, closed-form solutions are derived for two infinite sets of the involved parameters. These finite-sum solutions are written in terms of elementary functions that in general are quasi-polynomials. The coefficients of the expansion obey a three-term recurrence relation which in some particular cases may become two-term. We discuss the case when the recurrence relation involves two non-successive terms and show that the coefficients of the expansion are then calculated explicitly and the general solution of the Heun equation is constructed as a combination of several hypergeometric functions with quasi-polynomial pre-factors.
\\
\\ \textbf{MSC numbers: }33E30 Other functions coming from differential, difference and integral equations, 34B30 Special equations (Mathieu, Hill, Bessel, etc.), 30Bxx Series expansions
\\
\\ \textbf{Keywords:} Linear ordinary differential equation, Heun equation, special functions, series expansions, recurrence relations
\end{abstract}

\maketitle

\qquad The Heun differential equation \cite{Heun-1} and its four confluent forms \cite{Ronveaux-2}-\cite{Olver-6} have been subject to intensive investigations during the last two decades due to a large number of applications in mathematics and physics, particularly, in quantum mechanics (see, e.g., \cite{Ronveaux-2} and references therein). This is a second-order ordinary linear differential equation that has four regular singular points located
at $z=0, 1, a$ and $\infty.$
The equation presents direct generalization of the Gauss hypergeometric equation and includes as particular cases many other known equations of mathematical physics,
e.g., the Lam\'{e} and Mathieu equations.

\qquad The closed-form solutions of the Heun equations in terms of simpler functions are rare. Many of such solutions are derived via termination of infinite series solutions in terms of known special functions, in particular, hypergeometric ones (see, e.g., \cite{Svartholm-7}-\cite{Ishkhanyan-15}). In the present paper, we introduce a series solution of the general Heun equation in terms of the incomplete Beta functions \cite{Abramowitz-16}. The expansion is applicable if one of the characteristic exponents of the singularity at infinity is zero. We show that the constructed expansion generates closed-form solutions for two infinite sets of the involved parameters. These solutions in general are quasi-polynomials. The coefficients of the expansion in general obey a three-term recurrence relation. We discuss the cases when the recurrence relation becomes two-term.
In a particular case when the recurrence relation involves two non-successive terms
the coefficients of the expansion are calculated explicitly and the general solution of the Heun equation is written in quadratures or in terms of several hypergeometric functions.

\qquad The canonical form of the general Heun equation is
\begin{equation}\label{1}
u^{''}+\left(\frac{\gamma}{z}+\frac{\delta}{z-1}+\frac{\varepsilon}{z-a}\right)u^{'}+\frac{\alpha \beta z-q}{z(z-1)(z-a)}u=0,
\end{equation}
where the involved parameters satisfy the Fuchsian condition $1+\alpha+\beta=\gamma+\delta+\varepsilon.$

\qquad Consider an expansion of the solution of this equation in the following form:
\begin{equation}\label{2}
u=\sum_{n=0}^{\infty}a_{n}u_{n}=\sum_{n=0}^{\infty}a_{n}B_{z}(\gamma_{n}, \delta_{n}),
\end{equation}
where $B_{z}(\gamma_{n}, \delta_{n})$ is the incomplete Beta function satisfying the following equation \cite{Abramowitz-16}:
\begin{equation}\label{3}
u_{n}^{''}+\left(\frac{1-\gamma_{n}}{z}+\frac{1-\delta_{n}}{z-1}\right)u_{n}^{'}=0.
\end{equation}
Using the known relations between the Beta functions
\begin{equation}\label{4}
B_{z}(\gamma, \delta)=\frac{1}{\gamma}[z^{\gamma}(1-z)^{\delta}+(\gamma+\delta)B_{z}(\gamma+1, \delta)],
\end{equation}
\begin{equation}\label{5}
B_{z}(\gamma, \delta)=\frac{1}{\gamma+\delta-1}[(\gamma-1)B_{z}(\gamma-1, \delta)-z^{\gamma-1}(1-z)^{\delta}],
\end{equation}
we have
\begin{equation}\label{6}
(1-z)u_{n}^{'}=(\gamma_{n}-1)u_{n-1}-(\gamma_{n}-1+\delta_{n})u_{n},
\end{equation}
\begin{equation}\label{7}
z(1-z)u_{n}^{'}=\gamma_{n}u_{n}-(\gamma_{n}+\delta_{n})u_{n+1},
\end{equation}
where we assumed $\gamma_{n\pm1}\equiv \gamma_{n}\pm1.$

\qquad Substitution of equations \eqref{2} and \eqref{3} into Eq. \eqref{1} gives
\begin{equation}\label{8}
\sum_{n=0}^{\infty}a_{n}[(A_{n}z(1-z)+D_{n}(1-z)+(1-a)C_{n})u_{n}^{'}+(\alpha\beta z-q)u_{n}]=0,
\end{equation}
where
\begin{equation}\label{9}
\left\{
\begin{array}{rcl}
A_{n}=-(\gamma-1+\gamma_{n})-\varepsilon-C_{n} \\
D_{n}=a(\gamma-1+\gamma_{n})+(a-1)C_{n} \\
C_{n}=(\delta-1+\delta_{n})
\end{array}
\right.
\end{equation}
\\A further observation now is that the term $C_{n}u_{n}^{'}+\alpha\beta z u_{n}$  is not expressed as a linear combination of functions $u_n$ for any (nonzero) values of involved parameters
so that it must be cancelled. Thus, we can treat only the case $\alpha\beta=0.$ Since $\alpha$ and $\beta$ are the characteristic exponents of the regular singularity at infinity,
we thus see that the expansion is possible only if an exponent of $z=\infty$ is zero.

\qquad Putting now $C_{n}=0,$ we get that for all $n$ should hold
\begin{equation}\label{10}
\delta_{n}=1-\delta.
\end{equation}
Eq. \eqref{8} then reads
\begin{equation}\label{11}
\sum_{n=0}^{\infty}a_{n}[(A_{n} z(1-z)+D_{n}(1-z))u_{n}^{'}-qu_{n}]=0,
\end{equation}
and the corresponding recurrence relation for the coefficients of the series \eqref{2} is
\begin{equation}\label{12}
a_{n} a(\gamma-1+\gamma_{n})(\gamma_{n}-1)+$$
$$a_{n-1}[-a(\gamma-2+\gamma_{n})(\gamma_{n}+\delta_{n}-2)-((\gamma-2+\gamma_{n})+\varepsilon(\gamma_{n}-1)-q]+$$
$$a_{n-2}[(\gamma_{n}+\delta_{n}-2)((\gamma-3+\gamma_{n})+\varepsilon)]=0.
\end{equation}
\\From the left-hand side termination conditions $a_{-2}=a_{-1}=0$ and $a_{0}-\forall,$ we find that should be $\gamma_{n}=1+n$ or $1-\gamma+n.$
The case $\gamma_{n}=1+n$ does not work since then for $n=0$ Eq. \eqref{6} involves $u_{-1}=B(z, 0, \delta_{n})$. So, we get
\begin{equation}\label{13}
\gamma_{n}=1-\gamma+n,
\end{equation}
hence, finally,
\begin{equation}\label{14}
u=\sum_{n=0}^{\infty}a_{n}B_{z}(1-\gamma+n, 1-\delta).
\end{equation}

\qquad The coefficients of the three term recurrence relation given by equation \eqref{12}:
\begin{equation}\label{15}
R_{n}a_{n}+Q_{n-1}a_{n-1}+P_{n-2}a_{n-2}=0,
\end{equation}
finally read:
\begin{equation}\label{16}
R_{n}=an(n-\gamma),
\end{equation}
\begin{equation}\label{17}
Q_{n}=-an(n+1-\gamma-\delta)-(n+\varepsilon)(n+1-\gamma)-q,
\end{equation}
\begin{equation}\label{18}
P_{n}=(n+2-\gamma-\delta)(n+\varepsilon).
\end{equation}
The series \eqref{2} is terminated if $a_{N}\neq0$ and $a_{N+1}=a_{N+2}=0$ for some $N=1,2,3,...$. The second condition is fulfilled if
\begin{equation}\label{19}
\varepsilon=-N, \qquad N=1, 2, ...
\end{equation}
or if
\begin{equation}\label{20}
\gamma+\delta-2=+N, \qquad N=1, 2, ...
\end{equation}
while the equation $a_{N+1}=0$ imposes an additional restriction on the accessory parameter $q$.

\qquad The first choice complements the expansion of \cite{Ishkhanyan-12} extending it to the case $\alpha\beta=0$. Indeed, it is not difficult to show that the equation for the accessory parameter $q$ is exactly the one which is derived from that of \cite{Ishkhanyan-12} by putting in it $\alpha\beta=0$. Note that this equation always has a trivial root $q=0$ so the number of nontrivial solutions for $\varepsilon=-N$, in general, is $N$. The second case of termination of the series, that is the case given by Eq. \eqref{20}, defines a different solution to the Heun equation.
In this case the number of nonzero values of the accessory parameter $q$ for which the termination occurs is $N+1$.

\qquad Note that, according to the lemma of \cite{Ishkhanyan-15}, the finite-sum solutions are always written in terms of elementary functions. However, generally these are not polynomials but quasi-polynomials. Concerning the infinite series, we note that the coefficients $a_{n}$ are not calculated explicitly unless the recurrence relation \eqref{15} is reduced to a two-term one. The latter may involve successive or non-successive terms. Eq. \eqref{15} is reduced to a two-term one involving two successive terms in several cases, which, however, are rather simple such as, e.g., the trivial case $a=0$. Different is the case when the recurrence relation involves two non-successive terms. For this to occur for some set of the involved parameters,
the coefficient $Q_{n}$ should identically vanish for all $n$. This is the case if
\begin{equation}\label{21}
a=-1,\qquad \delta=-\varepsilon,\qquad q=(\gamma-1) \varepsilon.
\end{equation}
The expansion coefficients are then explicitly calculated in terms of the Gamma functions. Using the Pochhammer symbol, the result reads
\begin{equation}\label{22}
a_{n}=\frac{1+(-1)^{n}}{2}c_{n/2}, c_{k}=\frac{(\varepsilon/2)_{k}(1+(\varepsilon-\gamma)/2)_{k}}{k!(1-\gamma/2)_{k}}.
\end{equation}
This leads to the following general solution of the Heun equation for parameters \eqref{21} along with $\alpha\beta=0$ involving quadratures:
\begin{equation}\label{23}
H(-1,(\gamma-1)\varepsilon; \alpha\beta=0,\gamma,\delta=-\varepsilon; z)=$$
$$C_{1} \int_{0}^{z}(1-z)^{\varepsilon}z^{-\gamma} {\cdot_{2}}F_{1}\left(\frac{\varepsilon}{2}, 1+\frac{\varepsilon-\gamma}{2}; 1-\frac{\gamma}{2};z^{2}\right)dz+$$
$$C_{2}\left(-\frac{\gamma}{q}+\int_{0}^{z}(1-z)^{\varepsilon}{\cdot_2}F_{1}\left(1+\frac{\varepsilon}{2}, \frac{\varepsilon+\gamma}{2}; 1+\frac{\gamma}{2};z^{2}\right)dz\right),
\end{equation}
where $C_{1,2}$ are arbitrary constants. Note that here $u(0)=-C_{2}\gamma/q.$ Notably, the value adopted by the solution at $z=1$ is
also explicitly calculated in terms of the Gamma functions and a Clausen generalized hypergeometric function:
\begin{equation}\label{24-1}
u(1)=C_{1}\frac{\Gamma \left(1-\frac{\gamma}{2}\right)\Gamma \left(\frac{1}{2}+\frac{\varepsilon}{2}\right)}{\sqrt{\pi}(1-\gamma)\Gamma \left(1+\frac{\varepsilon}{2}-\frac{\gamma}{2}\right)}+$$
$$C_{2}\left(-\frac{\gamma}{q}+ {_3}F_{2}\left(\frac{1}{2}, 1, \frac{\varepsilon+\gamma}{2};1+\frac{\gamma}{2}, \frac{3}{2}+\frac{\varepsilon}{2}; 1\right)\right).
\end{equation}
The parametric excess for the involved  Clausen's hypergeometric function is $1$. Then if additionally an upper parameter is a negative integer the function becomes Saalschutzian, and it is expressed in terms of the Gamma functions \cite{Snow-5}. Hence, if $\gamma+\varepsilon=-2N, \ N=1, 2, ...,$  the value $u(1)$ is eventually calculated in terms of the Gamma functions.
Note that this condition differs from the conditions \eqref{19}-\eqref{20} for termination of the above series.

\qquad A further observation now is that the integrals involved in Eq. \eqref{23} can be explicitly calculated applying integration by parts. The final result is given by a rather cumbersome expression involving several hypergeometric functions of the argument $z^{2}$.
This can be straightforwardly understood if the solution \eqref{23} is rewritten, using Eq. \eqref{1}, as
\begin{equation}\label{24-2}
u=-\frac{z(z-1)(z-a)}{\alpha\beta z-q}\left(v'+\left(\frac{\gamma}{z}+\frac{\delta}{z-1}+\frac{\varepsilon}{z-a}\right)v\right),
\end{equation}
or, equivalently, for parameters \eqref{21} along with $\alpha\beta=0$:
\begin{equation}\label{25}
u=\left(z(z^{2}-1)v'+(\gamma(z^{2}-1)-2\varepsilon z)v\right)/q,
\end{equation}
where
\begin{equation}\label{26}
v=C_{1}(1-z)^{\varepsilon}z^{-\gamma}{\cdot_2}F_{1}\left(\frac{\varepsilon}{2}, 1+\frac{\varepsilon-\gamma}{2}; 1-\frac{\gamma}{2}; z^{2}\right)+$$
$$C_{2}(1-z)^{\varepsilon}{\cdot_2}F_{1}\left(1+\frac{\varepsilon}{2}, \frac{\varepsilon+\gamma}{2}; 1+\frac{\gamma}{2}; z^{2}\right).
\end{equation}
It is seen that solution \eqref{23} presents a combination of four hypergeometric functions (a pair of functions for each of the involved fundamental solutions) with quasi-polynomial pre-factors.

\qquad It is interesting to compare the derived solution with the known result for $a=-1,$ $\delta=\varepsilon$ and $q=0$
(instead of $\alpha\beta=0,$) which is also written in terms of the Gauss hypergeometric functions of the argument $z^{2}$ \cite{Maier-17}-\cite{Vidunas-18}:
\begin{equation}\label{27}
H(-1, 0; \alpha, \beta, \gamma, \delta=\varepsilon; z)=$$
$$C_{1}z^{1-\gamma}\cdot_ {2}F_{1}\left(\varepsilon-\frac{\alpha}{2}, \varepsilon-\frac{\beta}{2}; \frac{3-\gamma}{2}; z^{2}\right) $$
$$+C_{2} \cdot _{2}F_{1}\left(\frac{\alpha}{2}, \frac{\beta}{2}; \frac{1+\gamma}{2}; z^{2}\right).
\end{equation}

\qquad Thus, we have presented a series solution of the general Heun equation in terms of the incomplete Beta functions. The expansion applies if a characteristic exponent of the singularity at infinity is zero. We have shown that closed-form solutions can be derived for two infinite sets of the involved parameters. These finite-sum solutions are written in terms of elementary functions that in general are quasi-polynomials. The coefficients of the expansion in general obey a three-term recurrence relation. We have mentioned a particular case when the recurrence relation is reduced to a two-term one involving two non-successive terms. The coefficients of the expansion are then calculated explicitly and
the general solution of the Heun equation is written in quadratures or as a combination of hypergeometric functions.

\qquad A concluding remark concerning the constructed expansion is that it can be significantly extended, if combined with the Darboux transformation \cite{Darboux-20} of the Heun equation or if the properties of the derivatives of the Heun functions are applied \cite{Ishkhanyan-12}, \cite{Fiziev-21}-\cite{Ishkhanyan-24}, to construct expansions in terms of higher transcendental functions, e.g., the Goursat \cite{Shahnazaryan-23} and Appell \cite{Ishkhanyan-24} generalized hypergeometric functions. Using the relationship between the Heun class of second-order linear equations and the Painlev\'{e}
second-order nonlinear equations \cite{Slavyanov-25}, these expansions can also be adopted to discuss several nonlinear problems.

\subsection*{Acknowledgments}

This research has been conducted within the scope of the International Associated Laboratory (CNRS-France \& SCS-Armenia) IRMAS. The research has received funding from the European Union Seventh Framework Programme (FP7/2007-2013) under great agreement No. 295025, - IPERA. The work has been supported by the Armenian State Committee of Science (SCS Grant No. 13RB-052).


\end{document}